\documentclass{article}
\usepackage[dvips]{graphicx}
\usepackage{amssymb}
\usepackage{amsmath}

\usepackage{caption}
\usepackage{hyperref}
\usepackage{arcs}
\usepackage {xcolor} 
\usepackage{subfigure}
\usepackage{authblk}
\usepackage{calrsfs}

\newtheorem{theorem}{Theorem}

\usepackage{graphicx} % Required for inserting images

\newenvironment{AMS}{\small\bf 2020 AMS subject classification: }{} 

\title{Polynomial-free unisolvence of polyharmonic splines with odd exponent by random sampling}

\author{Alvise Sommariva, Marco Vianello} 
\affil{University of Padova, Italy}

\begin{document}

\maketitle

\begin{abstract}
In a recent paper almost sure unisolvence of RBF interpolation at random points with no polynomial addition was proved, for Thin-Plate Splines and Radial Powers with noninteger exponent. The proving technique left unsolved the case of odd exponents. In this short note   we prove almost sure polynomial-free unisolvence in such instances, by a deeper analysis of the interpolation matrix determinant and fundamental properties of analytic functions.   
\end{abstract}

\vskip0.2cm
\noindent
\begin{AMS}
{\rm 65D05,65D12.}
\end{AMS}
\vskip0.2cm
\noindent
{\small{\bf Keywords:} multivariate interpolation, Radial Basis Functions, polyharmonic splines, odd integer exponent, unisolvence, analytic functions.}

\section{Introduction}

Interpolation by Radial Basis Functions (RBF) is nowadays one of the basic tools of scattered data approximation and meshfree methods. In the case of Positive Definite RBF, such as Gaussians or Inverse Multiquadrics, unisolvence is guaranteed by the fact that the interpolation matrix is (symmetric) positive definite. Traditionally, the case of Conditionally Positive Definite (CPD) RBF of order $m$, such as MultiQuadrics and Polyharmonic Splines, is treated by adding a suitable polynomial term of degree $m-1$, ensuring that the interpolation matrix becomes  positive definite; cf., e.g., \cite{F07,I03}. For CPD-RBF of order $m=1$, such as MultiQuadrics and distance functions, it is however theoretically known that the polynomial-free interpolation matrix is nonsingular itself in any dimension, cf. \cite{M86}. Unisolvence was also proved in the very special case of univariate cubic powers of distance functions \cite{BS87}.

Though it has been known for a long time that Polyharmonic Splines can have nonsingular polynomial-free interpolation matrices in many applications, cf. e.g. \cite{P22} with the references therein,  their polynomial-free unisolvence has not been studied theoretically. It is worth recalling the following statement that appeared in the popular treatise \cite{F07}: {\em ``There is no result that states that interpolation with Thin-Plate Splines (or any other strictly conditionally positive definite function of order $m\geq 2$) without the addition of an appropriate degree $m-1$ polynomial is well-posed''}, and the situation did not apparently change until very recently. 

Indeed, polynomial-free unisolvence has been reconsidered experimentally in \cite{P22}. On the other hand, two meaningful theoretical results have been proved in the framework of random sampling. 
The first concerns interpolation by fixed centers \cite{DASV23}, whereas the second
treats the standard and more difficult case of centers coinciding with the sampling points \cite{BSV23}. In both cases one of the key aspects is that Polyharmonic Splines $\phi(\|x-\overline{x}\|_2)$, that correspond to the radial functions 
$$
\phi(r)=r^{2k} \log(r)\;,\;\;k\in \mathbb{N}\;\;\mbox{(TPS: Thin-Plate Splines, order $m=k+1$)}
$$
and 
$$
\phi(r)=r^\nu\;,\;\;0<\nu\notin 2\mathbb{N}\;\;\mbox{(RP, Radial Powers, order $m=\lceil\nu/2\rceil$)}\;,
$$ 
are {\em real analytic functions} off their center $\overline{x}$, due to analyticity of the univariate functions $\phi(\cdot)$ and $\sqrt{\cdot}$ in $\mathbb{R}^+$. Concerning the role of analyticity in unisolvence by random sampling it is also worth quoting the recent paper \cite{XN23}, where nevertheless only spaces of fixed analytic functions independent of the sampling points were considered. 

However, the proving technique in \cite{BSV23} was able to treat TPS and RP with noninteger exponent, but left unsolved the most usual case for RP, i.e. that of {\em odd integer exponents}. In this brief note we fill the gap, by a deeper analysis of the interpolation matrix determinant, resorting to some fundamental properties of analytic functions. The result is stated in the following:

\begin{theorem}
Let $\Omega$ be an open connected subset of $\mathbb{R}^d$, $d\geq 2$, and 
$\{x_i\}_{i\geq 1}$ be a randomly distributed sequence on $\Omega$ with respect to any given probability density $\sigma(x)$, i.e. a point sequence produced by sampling a sequence of absolutely continuous random variables $\{X_i\}_{i\geq 1}$ which are independent and identically distributed in $\Omega$ with density $\sigma\in L^ 1_+(\Omega)$. Moreover, let $V_n=\left[\phi(\|x_i-x_j\|_2)\right]$, $1\leq i,j\leq n$, $\phi(r)=r^\nu$, be the interpolation matrix with respect to RP with odd integer exponent $\nu=2k+1$, $0\leq k\in \mathbb{N}$. 

Then, for every $n\geq 2$ the matrix $V_n$ is {\em a.s.}\ (almost surely) nonsingular.   
\end{theorem}
\vskip0.5cm 
\noindent
{\bf Proof.} Let us define for convenience $\phi_j(x)=\phi(\|x-x_j\|_2)$, and 
$$
V(\phi_1,\dots,\phi_n;y_1,\dots,y_n)=[\phi_j(y_i)]\;,\;\;1\leq i,j \leq n\;,
$$
so that $V_n=V(\phi_1,\dots,\phi_n;x_1,\dots,x_n)$. Notice that the diagonal of $V_n$ is zero, since $\phi(0)=0$. Now, the functions $\{\phi_j(x)\}$ are linearly independent in $\Omega$ if and only if the points $\{x_j\}$ are distinct. Clearly, this is a necessary condition for unisolvence. In fact, if the functions $\{\phi_j(x)\}$ were linearly dependent, one of them would be linear combination of the others (which are analytic in $\Omega$ off their centers) and thus would become analytic at its own center.

It is also worth recalling that a subset of $\Omega$ has null measure with respect to $d\sigma=\sigma(x)\,dx$, if it has null Lebesgue measure (i.e., it is a so-called ``null set''). 

First, we prove by induction on $n\geq 2$ that 
\begin{itemize}
\item [$(i)$] 
$x_1,\dots,x_n$ are a.s. distinct;
\item [$(ii)$] 
the subdeterminant $\delta_{n-1}=det(V(\phi_1,\dots,\phi_{n-1};x_2,\dots,x_n))$ is a.s. nonzero. 
\end{itemize}
Then, by $(i)$ and $(ii)$ we shall prove that $V_n$ itself is a.s. nonsingular. 

The statements $(i)$ and $(ii)$ hold for $n=2$. Indeed, the probability that $x_2=x_1$ given $x_1$ is zero, since the probability that a random point belongs to any finite set is zero (a finite set being a null set). This entails that  $\delta_1=det(V(\phi_1;x_2))=\phi_1(x_2)$ is a.s. positive. 

We make now the inductive step. The points $\{x_1,\dots,x_{n+1}\}$ are a.s. distinct because such are $\{x_1,\dots,x_n\}$ and the probability that $x_{n+1}$ coincides with one of them is zero, a finite set being a null set. That is, $(i)$ holds for $n+1$. 

As for $(ii)$, consider the matrix  
$$
W(x)=\left(\begin{array} {cccccc}
\phi_1(x_2) & 0 & \phi_3(x_2) & \cdots & \phi_{n-1}(x_2) & \phi_n(x_2)\\
\\
\phi_1(x_3) & \phi_2(x_3) & 0 & \cdots & \phi_{n-1}(x_3) & \phi_n(x_3)\\
\\
\vdots & \vdots & \vdots & \ddots & \vdots  & \vdots\\
\\
\phi_1(x_{n-1}) & \phi_2(x_{n-1}) & \phi_3(x_{n-1}) & \cdots & 0  & \phi_n(x_{n-1})\\
\\
\phi_1(x_{n}) & \phi_2(x_{n}) & \phi_3(x_{n}) & \cdots & \phi_{n-1}(x_n) & 0\\
\\
\phi_1(x) & \phi_2(x) & \phi_3(x) & \cdots & \phi_{n-1}(x) & \phi_n(x)\\
\\
\end{array} \right)
$$
and observe that $W(x_{n+1})=V(\phi_1,\dots,\phi_{n-1},\phi_n;x_2,\dots,x_n,x_{n+1})$. Developing the determinant of $W(x)$ by Laplace rule along the last row we get $$G(x)=det(W(x))=\delta_{n-1}\phi_n(x)+c_{n-1}\phi_{n-1}(x)+\dots+c_1\phi_1(x)\;,$$ where the $\{c_j\}$ are the corresponding minors with the appropriate sign, that do not depend on $x$. Therefore, $G$ is an analytic function in the open connected set $\Omega\setminus \{x_1,\dots,x_n\}$, being $G\in span(\phi_1,\dots,\phi_n)$. Moreover, $G$ is a.s. not identically zero in $\Omega$, because $\delta_{n-1}$ is nonzero and the $\{\phi_j\}$ are linearly independent in $\Omega$, by inductive hypothesis. On the other hand, by continuity in $\Omega$, $G(x)$ is a.s. not identically zero also in $\Omega\setminus \{x_1,\dots,x_n\}$. 

Then, 
$\delta_n=det(U(x_{n+1}))=G(x_{n+1})$ is a.s. nonzero, 
since the zero set of a not identically zero real analytic function on an open connected set in $\mathbb{R}^d$ is a null set (cf. \cite{M20} for an elementary proof). More precisely, 
denoting by $Z_G$ the zero set of $G$ in $\Omega$, 
we have that $$Z_G=(Z_G\cap \{x_1,\dots,x_n\})\cup 
(Z_G\cap (\Omega\setminus \{x_1,\dots,x_n\}))\;.$$ 
Hence $Z_G$ is a null set if $G\not\equiv 0$, because the first intersection is a finite set, and the second is the zero set of a not identically zero real analytic function. Considering the probability of the corresponding events we can then write 
$$
prob\{G(x_{n+1})=0\}=prob\{G\equiv 0\}+prob\{G\not\equiv 0\;\&\;x_{n+1}\in Z_G\}
=0+0=0\;.
$$ 
That is, also $(ii)$ holds for $n+1$ and the inductive step is completed.

We can now prove that $V_n$ is a.s. nonsingular for $n\geq 2$. The assertion is clearly true for $n=2$, since $det(V_2)=\phi_2(x_1)\phi_1(x_2)=-\phi_1^2(x_2)\neq 0$. For $n\geq 3$, consider the $n\times n$ matrix 
$$
U(x)=\left(\begin{array} {cccccc}
0 & \phi_2(x_1)  & \cdots 
& \phi_{n-2}(x_1) & \phi_{n-1}(x_1)  & \phi_1(x)\\
\\
\phi_1(x_2) & 0  & \cdots & \phi_{n-2}(x_2) & \phi_{n-1}(x_2)  & \phi_2(x)\\
\\
\vdots & \vdots  & \ddots & \vdots  & \vdots & \vdots\\
\\
\phi_1(x_{n-2}) & \phi_2(x_{n-2})  & \cdots & 0 & \phi_{n-1}(x_{n-2})  & \phi_{n-2}(x)\\
\\
\phi_1(x_{n-1}) & \phi_2(x_{n-1})  & \cdots & \phi_{n-2}(x_{n-1}) & 0  & \phi_{n-1}(x)\\
\\
\phi_1(x) & \phi_2(x)  & \cdots & \phi_{n-2}(x) & \phi_{n-1}(x) & 0\\
\\
\end{array} \right)
$$

Applying Laplace determinantal rule to the $n$-th row, we get 
$$
F(x)=det(U(x))=\alpha_1(x)\phi_1(x)+\dots
+\alpha_{n-1}(x)\phi_{n-1}(x)\;,
$$
where $\alpha_1,\dots,\alpha_{n-1}$ are the corresponding minors with the appropriate sign and clearly $\alpha_j\in span(\phi_1,\dots,\phi_{n-1})$. 
The claim is that $F(x)$ is a.s. not identically zero in $\Omega$. 

Expanding the computation of $\alpha_1(x)$ and $\alpha_{n-1}(x)$, 
it is easy to see that 
$$
F(x)=-d_{n-2}\phi_{n-1}^2(x)+A(x)\phi_{n-1}(x)+B(x)\;,\;\;d_{n-2}=det(V_{n-2})\;,
$$
where $A\in span(\phi_1,\dots,\phi_{n-2})$ and $B(x)\in span(\phi_j\phi_k\,,\,1\leq j,k\leq n-2)$. More precisely, developing $\alpha_1(x)$ and $\alpha_{n-1}(x)$ 
by Laplace rule using the last column we obtain  
$$
A(x)=(-1)^{n+1}\phi_1(x)\,det(V(\phi_2,\dots,\phi_{n-1};x_1,\dots,x_{n-2}))
$$
$$
+(-1)^{3n-1}\phi_1(x)\,det(V(\phi_1,\dots,\phi_{n-2};x_2,\dots,x_{n-1}))+C(x)
$$
where $C\in span(\phi_2,\dots,\phi_{n-2})$, from which follows
$$
A(x)={(-1)}^{n+1}2\delta_{n-2}\phi_1(x)+C(x)\;,
$$
since $V(\phi_2,\dots,\phi_{n-1};x_1,\dots,x_{n-2})=V^t(\phi_1,\dots,\phi_{n-2};x_2,\dots,x_{n-1})$ and hence the two matrices above have the same determinant. 

Notice that $A$ is a.s. not identically zero in $\Omega$, because by inductive hypothesis $\delta_{n-2}$ is nonzero and $\phi_1,\dots,\phi_{n-1}$ are linearly independent. 
And more, being continuous in $\Omega$ and analytic in the open connected set $\Omega\setminus \{x_1,\dots,x_{n-2}\}$, it is a.s. not identically zero also there, otherwise by continuity it would be a.s. identically zero on the whole $\Omega$. Then, it is a.s. not identically zero in a neighborhood of $x_{n-1}$, because its zero set must be a null set in $\Omega\setminus \{x_1,\dots,x_{n-2}\}$ (cf. e.g. \cite{M20}).

Assume now that $F\equiv 0$. Then we would have $d_{n-2}\phi_{n-1}^2-B\equiv A\phi_{n-1}$. This leads to a contradiction, since $d_{n-2}\phi_{n-1}^2$ is a polynomial of degree $2\nu$ 
and $B$ is analytic in a neighbourhood of $x_{n-1}$, so that 
$d_{n-2}\phi_{n-1}^2-B$ is analytic in such a neighbourhood, whereas $A\phi_{n-1}$ has a.s. a singularity at $x_{n-1}$. 

To prove the latter assertion, by $A\not\equiv 0$ it follows that, in the direction of some unit vector $u$, the univariate analytic function $\alpha(t)=A(x_{n-1}+tu)$ is not identically zero in a neighbourhood of $t=0$. On the other hand, $\phi_{n-1}(x_{n-1}+tu)=\phi(t)=|t|^\nu$ has a discontinuity of the $\nu$-th derivative.

If $\alpha(0)\neq 0$, by Leibniz rule for the derivatives of a product this leads immediately to the fact that $\alpha\phi$ has a discontinuity of the $\nu$-th derivative at $t=0$. If $\alpha(0)=0$, by a well-known result on the zeros of real analytic functions, $\alpha(t)=t^k\beta(t)$, where $k$ is called the order of the zero, $1\leq k\in \mathbb{N}$, and $\beta$ is (locally) analytic with $\beta(0)\neq 0$ (cf., e.g., \cite{KP02}). Then, again by Leibniz rule, $\alpha\phi$ has a discontinuity of the $(\nu+k$)-th derivative at $t=0$. In any case, $A\phi_{n-1}$ has a.s. a singularity at $x_{n-1}$.

Observe that $U(x_{n})=V_{n}$ since $\phi_{n}(x_{n})=0$ and $\phi_j(x_{n})=\phi_{n}(x_{j})$ for $j=1,\ldots,n-1$.  
On the other hand, $F$ is continuous in $\Omega$ and analytic in the open connected set $\Omega\setminus \{x_1,\dots,x_{n-1}\}$, being sum of products of analytic functions, and is a.s. not identically zero also there, otherwise by continuity it would be a.s. identically zero on the whole $\Omega$. 

Then, defining $Z_F$ as above for $G$, we can conclude in the same way that $prob\{F(x_{n})=0\}=0$ and thus $det(V_n)=F(x_{n})$ is a.s. nonzero.\hspace{0.2cm} $\square$

\section*{Acknowledgements}

Work partially
supported by the
DOR funds 
of the University of Padova, and by the INdAM-GNCS.
This research has been accomplished within the RITA ``Research ITalian network on Approximation", the SIMAI Activity Group ANA\&A and the UMI Group TAA ``Approximation Theory and Applications".


\begin{thebibliography}{99}

\bibitem{BS87} L.P. Bos, K. Salkauskas, On the matrix $[|x_i-x_j|^3]$ and the cubic spline continuity equations, J. Approx. Theory 51 (1987), 81--88. 

\bibitem{BSV23} L.P. Bos, A. Sommariva, M. Vianello, A note on polynomial-free unisolvence of polyharmonic splines at random points, preprint available at  \url{https://arxiv.org/abs/2312.13710}.

\bibitem{DASV23} F. Dell'Accio, A. Sommariva, M. Vianello, 
Random sampling and unisolvent interpolation by almost everywhere analytic functions, Appl. Math. Lett. 145 (2023), 108734. 

\bibitem{F07} G.E. Fasshauer, Meshfree Approximation Methods with Matlab, Interdisciplinary Mathematical
Sciences, Vol. 6, World Scientific, 2007.

\bibitem{I03} A. Iske, On the Approximation Order and Numerical Stability of Local Lagrange Interpolation by Polyharmonic Splines, in: Modern Developments in Multivariate Approximation, Haussmann, W., Jetter, K., Reimer, M., Stöckler, J. (eds), Birkh{\"a}user Basel, 2003, 153--165.

\bibitem{KP02} S.G. Krantz and H.R. Parks, A Primer of Real Analytic Functions, Second Edition, Birkh\"auser, Boston, 2002.

\bibitem{M86} C.A. Micchelli, Interpolation of scattered data: distance matrices and conditionally positive definite functions, Constr. Approx. 2 (1986), 11--22.

\bibitem{M20} B.S. Mityagin, The Zero Set of a Real Analytic Function, Math. Notes 107 (2020), 529--530. 

\bibitem{P22} A. Pasioti, On the Constrained Solution of RBF Surface
Approximation, MDPI Mathematics 10 (15) (2022), 2582. 

\bibitem{XN23} Y. Xu, A. Narayan, Randomized weakly admissible meshes, 
J. Approx. Theory 285 (2023), 105835.

\end{thebibliography}
\end{document}